\documentclass[12pt]{article}
\usepackage{epsfig}

\title{Fast Computation of the  Expected  Loss  of a Loan Portfolio Tranche in the Gaussian Factor Model: Using Hermite Expansions for Higher Accuracy.}

\author{ Pavel Okunev\footnote{This work was supported by the Director, Office of Science, Office  of Advanced Scientific Computing Research, of the U.S. Department of Energy under Contract No. DE-AC03-76SF00098.} \footnote{E-mail: pokunev@math.lbl.gov } \\  Department of Mathematics\\LBNL and UC Berkeley\\ Berkeley, CA 94720}

\date{June 19,2005}

\begin{document}
\maketitle

\begin{abstract}
We propose a fast algorithm for computing the expected tranche loss in the Gaussian  factor model. We test it on portfolios ranging in size from 25 ( the size of DJ iTraxx Australia) to  100 (the size of DJCDX.NA.HY) with a single factor Gaussian model and show that the algorithm gives accurate results.  The algorithm proposed here is an extension of the algorithm proposed in \cite{PO}. The advantage of the new algorithm is that it works well for portfolios of smaller size for which the normal approximation proposed in \cite{PO} in not sufficiently accurate. The algorithm  is intended as an alternative to the much slower Fourier transform based methods \cite{MD}.
\end{abstract}

\section{The Gaussian Factor Model}

Let us consider a portfolio of $N$ loans. Let the notional of loan $i$ be equal to the fraction $f_i$ of the notional of the whole portfolio. This means that if loan $i$ defaults  and the entire notional of the loan is lost the portfolio loses fraction $f_i$ or $100f_i\%$  of its value. In practice when a loan $i$ defaults a fraction $r_i$ of its notional will be recovered by the creditors.  Thus the actual loss given default (LGD) of loan $i$ is
\begin{equation}
LGD_i=f_i(1-r_i)
\end{equation}
fraction or
\begin{equation}
LGD_i=100f_i(1-r_i)\%
\end{equation}
of the notional of the entire portfolio.

We now describe the Gaussian m-factor model of portfolio losses from default. The model requires a number of input parameters. For each loan $i$ we are give a probability $p_i$ of its default. Also for each  $i$ and each $k=1,\ldots,m$ we are given a number $w_{i,k}$ such that $\sum_{k=1}^m w_{i,k}^2<1$. The number  $w_{i,k}$  is the loading factor of the loan $i$ with respect to factor $k$. Let $\phi_1, \ldots, \phi_m$ and $\phi^i, i=1,\ldots,N$ be independent standard normal random variables. Let $\Phi(x)$ be the cdf of the standard normal distribution.
In our model loan $i$ defaults if 
\begin{equation}
\sum_{k=1}^m w_{i,k}\phi_k+\sqrt{1-\sum_{k=1}^m w_{i,k}^2}\phi^i<\Phi^{-1}(p_i)
\end{equation}
This indeed happens with probability $p_i$.
  The factors $\phi_1,\ldots,\phi_m$ are usually interpreted  as the state of the global economy, the state of the  regional economy, the state of a particular industry and so on. Thus they are  the factors that affect the default behavior of all or at least a large group of loans in the portfolio. The factors $\phi^1,\ldots,\phi^N$ are interpreted as the idiosyncratic risks of the loans in the portfolio.
 
Let $I_i$ be defined by
\begin{equation}
I_i=I_{\{loan \ i\ defaulted\}}
\end{equation}
  We define the random loss caused by the default of loan $i$ as 
\begin{equation}
L_i=f_i(1-r_i)I_i,
\end{equation}
where $r_i$ is the recovery rate of loan $i$.
The total loss of the portfolio is 
\begin{equation}
L= \sum_i L_i
\end{equation}

An important property of the Gaussian factor model is that tthe $L_i$'s are not independent of each other. Their mutual  dependence is induced by the dependence of each $L_i$ on the common factors $\phi_1,\ldots,\phi_m$. Historical data supports the conclusion that losses due to defaults on different loans are correlated with each other. Historical data can also be used to calibrate the loadings $w_{i,k}$.
 
\section{Conditional Portfolio Loss $L$}

When the values of the factors $\phi_1,\ldots,\phi_m$ are fixed, the probability of the default of loan $i$ becomes
\begin{equation}
p^i=\Phi^{-1}\left( \frac{p_i-\sum_kw_{i,k}\phi_k}{\sqrt{1-\sum_kw_{i,k}^2}} \right)
\end{equation}

The random losses $L_i$ become conditionally independent Bernoulli variables with the mean  given by 
\begin{equation}
E_{cond}(L_i)=f_i(1-r_i)p^i
\end{equation}
and the variance given by
\begin{equation}
VAR_{cond}(L_i)=f_i^2 (1-r_i)^2p^i(1-p^i)
\end{equation}   

By the Central Limit Theorem the conditional distribution of the portfolio loss $L$, given the values of the factors $\phi_1,\ldots,\phi_m$,  can be approximated by the normal distribution with the mean
\begin{equation}
\label{m}
E_{cond}(L)=\sum_i E_{cond}(L_i)
\end{equation}
and the variance

\begin{equation}
\label{var}
VAR_{cond}(L)=\sum_i VAR_{cond}(L_i)
\end{equation} 

In \cite{PO} it was shown that for portfolios of  125 names this approximation leads to accurate results. 

If the size of the portfolio is smaller than 125, for example 30 ( the size of DJ iTraxx ex Japan) or 50 (the size of DJ iTraxx CJ), then the Central Limit Theorem no longer provides a sufficiently accurate approximation to the conditional distribution of the portfolio loss $L$.
An  accurate representation of the  conditional distribution of the portfolio loss $L$ is given by its Hermite series expansion. For historical reasons this expansion is also known as the Charlier series expansion \cite{FE}, \cite{CA}.

\section{The Hermite Expansion of the Conditional Distribution of the Portfolio Loss $L$}

Let $F(x)$ be the c.d.f. of the  conditional distribution of the portfolio loss $L$. So that 
\begin{equation}
P(L\leq x)=F(x)
\end{equation}
 For  each fixed value of the factors $\phi_1,\ldots,\phi_m$ we define the normalized conditional loss $\tilde{L}$ by
\begin{equation} 
\tilde{L}=\frac{L-E_{cond}(L)}{\sqrt{VAR_{cond}(L)}}
\end{equation}
Let $\tilde{F}(x)$ be the c.d.f. of the  distribution of the normalized   conditional portfolio loss $\tilde{L}$. So that 
\begin{equation}
P(\tilde{L}\leq x)=\tilde{F}(x)
\end{equation}
We define the Hermite polynomial $H_n(x)$ of degree $n$ by
\begin{equation}
H_n(x)=(-1)^ne^{\frac{x^2}{2}}\frac{d^n}{dx^n}e^{\frac{-x^2}{2}}
\end{equation}
Let $c_n$ be defined by
\begin{equation}
c_n=\frac{(-1)^n}{n!} \int_{-\infty}^{\infty}H_n(x)d\tilde{F}(x)
\end{equation}
Then we have 
\begin{equation}
\label{exp}
\tilde{F}(x)=\sum_{i=0}^{\infty} \int_{-\infty}^{x}c_i H_i(t)\frac{e^{\frac{-t^2}{2}}}{\sqrt{2\pi}}dt
\end{equation}
The series above converges in the sense of distributions (generalized functions) \cite{R}. A good reference on the theory of distributions (generalized functions) is \cite{R}.
Let us pick a finite $N$. Then we have
\begin{equation}
\label{expN}
\tilde{F}(x) \approx \sum_{i=0}^{N} c_i \int_{-\infty}^{x}H_i(t)\frac{e^{\frac{-t^2}{2}}}{\sqrt{2\pi}}dt
\end{equation}
As before the approximation is in the sense of  generalized functions. 
Equation (\ref{expN}) implies that the distribution of the normalized   conditional portfolio loss $\tilde{L}$
can be approximated by a distribution with the density
\begin{equation}
\tilde{\rho}(x)=\sum_{i=0}^{N} c_i H_i(x)\frac{e^{\frac{-x^2}{2}}}{\sqrt{2\pi}}
\end{equation}
The function $\tilde{\rho}(x)$ is not necessarily nonnegative and therefore may not be a probability density in the strict sense. However, as is  explained in \cite{R}, this does not affect the validity of our final result (\ref{theformula}). Therefore we may treat $\tilde{\rho}(x)$ as a real probability density.

The distribution of the unnormalized loss $L$ can be approximated by a distribution with density 
\begin{equation}
\label{den1}
\rho(x)=\sum_{i=0}^{N} \frac{c_i}{\sqrt{VAR_{cond}(L)}} H_i\left(\frac{x-E_{cond}(L)}{\sqrt{VAR_{cond}(L)}}\right)\frac{e^{\frac{-\left(\frac{x-E_{cond}(L)}{\sqrt{VAR_{cond}(L)}}\right)^2}{2}}}{\sqrt{2\pi}}
\end{equation}

The joint distribution of the factors $\phi_1,\ldots,\phi_m$ and the portfolio loss $L$ can be approximated by a distribution with density
\begin{equation}
\label{den}
\rho_{joint}(\phi_1,\ldots,\phi_m,L)=\rho(L)\prod_{k=1}^m\rho_{G,0,1}(\phi_k),
\end{equation}
where $\rho_{G,0,1}(x)$ stands for the Gaussian density with mean $0$ and variance $1$.

Observe that the coefficient $c_n$  depends only on the moments of the distribution $\tilde{F}(x)$.
Since $L_i$'s are independent Bernoulli random variables these moments  are known analytically. Thus in the case under consideration all the $c_n$'s are known analytically.

If in equation (\ref{den1}) we set $N=1$ we  obtaine the standard approximation by the normal density proposed in \cite{PO}. Thus the algorithm proposed here is a generalization of the algorithm in \cite{PO}. We show later that it gives good numerical results even when the portfolio size is too small for the normal approximation to be accurate.

\section{Expected Loss of a Tranche of Loan Portfolio}

Let $0\leq a<b\leq 1$. We define a tranche loss profile $Tl_{a,b}(x)$ by
\begin{equation}
Tl_{a,b}(x)=\frac{min(b-a,max(x-a,0))}{b-a}
\end{equation}
Number $a$ is called the attachment point of a tranche, while $b$ is called the detachment point  of a  tranche. 
The expected loss of a tranche is then
\begin{equation}
TLoss(a,b)=\int Tl_{a,b}(L)\rho_{joint}(\phi_1,\ldots,\phi_m,L)d\phi_1\ldots\phi_mL
\end{equation}
This can be rewritten as a double integral
\begin{equation}
\label{theformula}
TLoss(a,b)=\int \int Tl_{a,b}(L)\rho(L)dL \prod_{k=1}^m\rho_{G,0,1}(\phi_k)d\phi_1\ldots\phi_m
\end{equation}
The inside integral with respect to $L$ can be  done analytically  for fixed values of the factors $\phi_1,\ldots,\phi_m$. The outside integral has to be computed numerically. However, since it is an integral of a bounded smooth function with respect to m-dimensional Gaussian density, it is one of the simpler integrals to compute numerically.

\section{Numerical Example}
In this section we apply the proposed algorithm to the single factor Gaussian model of a portfolio with $n$ names. We take $n$ to be 25 (size of DJ iTraxx Australia), 30  (size of DJ iTraxx ex Japan), 50 (size of DJ iTraxx CJ) and 100 (size of DJCDX.NA.HY).  We choose a single factor model because it is the one most frequently used in practice.  For each $n$ we compute the loss of the equity tranche  with the attachment point $a=0$ or $a=0\%$ and the detachment point 3\%.  The parameters of the porfolio are
\begin{eqnarray}
f_i&=&\frac{1}{n} \nonumber \\
p_i&=&0.015+\frac{0.05(i-1)}{n-1} \nonumber \\
r_i&=&0.5-\frac{0.1(i-1)}{n-1} \nonumber \\
w_{i1}&=&0.5-\frac{0.1(i-1)}{n-1},
\end{eqnarray}
where $i=1,\ldots,n$. Finally, we choose $N=5$ in (\ref{expN}).

\begin{figure}[htbp]
\caption{Equity Tranche Loss in the Gaussian Single Factor Model}
\epsfig{file=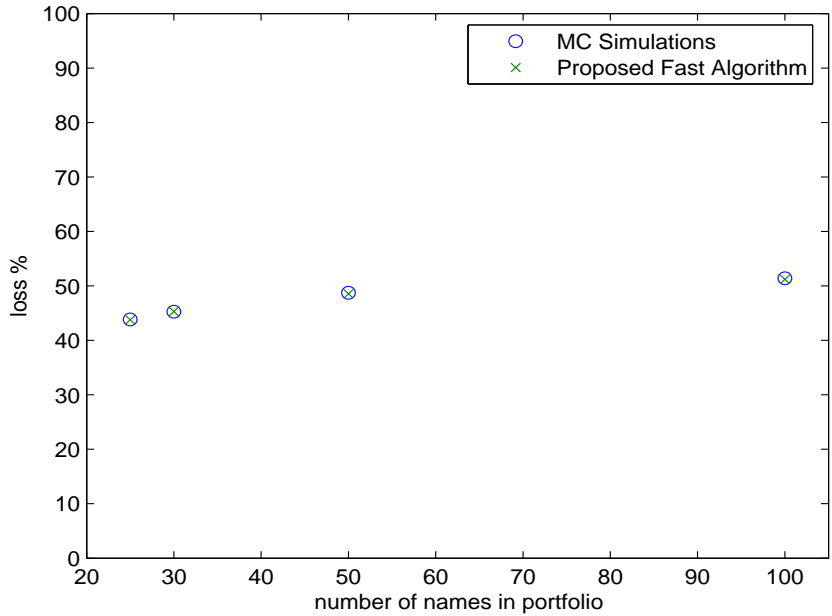,width=5.0 in,height=3.5 in}
\label{resultsfig}
\end{figure}

In Figure \ref{resultsfig} we compare the expected loss computed using $10^6$ Monte Carlo samples with the expected loss computed using formula (\ref{theformula}).\footnote{The author has the code implementing the algorithm described here in MATLAB, VBA for Excel and C.} The agreement between the two is good.

\section{Conclusions.}

To obtain the results in Figure \ref{resultsfig} we only needed to perform a single one dimensional numerical integration for each tranche. This is an  improvement over the Fourier transform based methods \cite{MD} which requires computing a large number of Fourier transforms for each tranche.  Each individual Fourier transform is as computationally expensive as (\ref{theformula}).

The expansion (\ref{expN}) is accurate even when the portfolio size is too small for the normal approximation of \cite{PO} to be precise. Thus we developed an algorithm which is as fast as the algorithm proposed in \cite{PO} but allows us to obtain higher precision for a portfolio of a given size  by including more terms in (\ref{expN}).
\section{Acknowledgments}
I thank my adviser A. Chorin for his help and guidance during my time in UC Berkeley. I thank Mathilda Regan and Valerie Heatlie for their help in preparing this article. I am also grateful to Ting Lei, Sunita Ganapati and George Wick  for encouraging my interest in finance.  Last, but not least, I thank my family for their constant support.


\begin{thebibliography}{99}

\bibitem{CA} H. Cramer. Mathematical Methods of Statistics.  Princeton University Press, 1954.

\bibitem{MD} A. Debuysscher, M. Szeg\"{o}, M. Freydefront and H. Tabe. Fourier Transform Method-Technical Document. Available from Moody's.

\bibitem{FE} W. Feller. An Introduction to Probability Theory and Its  Applications. Wiley, 1968.

\bibitem{PO} P.Okunev. A Fast Algorithm for Computing Expected Loan Portfolio Tranche Loss in the Gaussian Factor Model. LBNL-57676.

\bibitem{R} M.Reed and B.Simon. Methods of Modern Mathematical Physics, v. 2. Academic press,1972.

\end{thebibliography}
\end{document}